\documentclass[letterpaper, 10 pt, conference]{ieeeconf}  

\usepackage{dsfont}
\usepackage{amssymb}
\usepackage{amsmath}
\usepackage{graphicx}
\usepackage{wrapfig}
\usepackage{color}
\usepackage{soul} 
\usepackage{cite}

\usepackage{times}

\newtheorem{remark}{Remark}

\def\<{\langle}
\def\>{\rangle}

\def \t {\theta}
\def \a {\alpha}

\def \T {\Theta}

\def \d {\delta}

\def \O {\Phi}

\def \O {\Phi}

\IEEEoverridecommandlockouts                              
\title{\LARGE \bf
Disentangling Drift- and Control- Vector Fields for Interpretable Inference of Control-affine Systems 
}
\author{Vignesh Narayanan$^{1}$, Wei Miao$^{1}$ and Jr-Shin Li$^{1}$
\thanks{*This work was supported in part by the NSF under the awards CMMI-1763070 and CMMI-1737818 and by the NIH under the award R01GM131403.}
\thanks{$^{1}$The authors are with the Department of Electrical and Systems Engineering, Washington University in St. Louis, St. Louis, MO, USA {\tt\small vignesh.narayanan, weimiao, jsli @wustl.edu.}}%
}
\begin{document}
\maketitle
\thispagestyle{empty}
\pagestyle{empty}
\begin{abstract}
Many engineered as well as naturally occurring dynamical systems do not have an accurate mathematical model to describe their dynamic behavior. However, in many applications, it is possible to probe the system with external inputs and measure the process variables, resulting in abundant data repositories. Using the time-series data to infer a mathematical model that describes the underlying dynamical process is an important and challenging problem. In this work, we propose a model reconstruction procedure for inferring the dynamics of a class of nonlinear systems governed by an input affine structure. In particular, we propose a data generation and learning strategy to decouple the reconstruction problem associated with the drift- and control- vector fields, and enable quantification of their respective contributions to the dynamics of the system. This learning procedure leads to an interpretable and reliable model inference approach. We present several numerical examples to demonstrate the efficacy and flexibility of the proposed method.
\end{abstract}
\section{Introduction}
Engineered systems encountered in many scientific domains are increasingly complex and highly nonlinear. In recent times, the processes and mechanisms previously confined to biology, social science, etc., are viewed through the prism of systems theory, and, as a consequence, many high-dimensional complex systems are emerging (\cite{ronquist2017algorithm,altafini2018signed}). A common challenge in both the engineered and natural systems is the lack of accurate mathematical models describing their dynamic behavior. The ability to model such complex systems is essential for understanding the underlying dynamical processes, and, thereby, enabling safe and efficient utilization of such systems in safety-critical application domains such as medical diagnosis and prognosis. 

From the early transfer-function models to the widely used state-space representations \cite{kalman_63_mathematical}, an enormous body of systems representation and modeling approaches using first principles has been developed and reported \cite{ljung1999system}. However, in many emerging applications precise models that describe the dynamics of the system are not available. Due to the advances in actuator and sensing technologies, in many of these applications, it is possible to externally perturb the system and record the evolution of the system states \cite{ronquist2017algorithm,JCNS_2019biophysically,proctor_16_dmdc}. Therefore, it is desired to reconstruct an interpretable and accurate model for these systems by using the time-series data. 

Inferring the system models has been an active research topic in the past several decades \cite{ioannou1982asymptotic,ks1990identification,kreisselmeier1990stabilized,ljung1999system}. Plenitude of results have been reported on system identification and model learning  for linear dynamical systems \cite{ljung1999system,proctor_16_dmdc} and nonlinear dynamical systems \cite{schoukens2019nonlinear,brunton_16_sindy} with or without external controls. For instance, in the traditional systems theory, inferring a dynamical system model that encodes and reproduces the relationship in a given input-output dataset has been studied as a \emph{realization} problem \cite{isidori1995nonlinear,brockett1976volterra}. Moreover, methods such as dynamic mode-decomposition (DMD), DMDc \cite{proctor_16_dmdc}, SINDy \cite{brunton_16_sindy} and adaptive identifiers \cite{ioannou1982asymptotic, sastry2011adaptive,ks1990identification} have been proposed to infer the dynamical equations using simulated or measurement time-series data, and several successful applications of these approaches have been reported (see for example \cite{schmid_11_dmd_application,schoukens2019nonlinear}). 

In general, existing methods for model reconstruction aim to find a set of mathematical equations that best fits the given data samples. However, in many applications, the data-samples obtained need not be rich enough to fully recover the dynamics of the system, and the experimental costs to generate data can impede this process. In this context, generating data-samples efficiently is an important task. For instance, in adaptive control, the concepts of \emph{sufficient richness} and \emph{persistency of excitation} are used to describe the type of perturbation signals and regression functions that are desired for parameter convergence in a model estimation and learning problem \cite{sastry2011adaptive,ks1990identification}. Similarly, in canonical reinforcement learning algorithms, this idea is encapsulated by \emph{exploration} \cite{sutton1998introduction}. Furthermore, learning a model for the controlled system is typically done by combining the unknown parameters of the drift and control vector fields, resulting in an abstract dynamic model that explains the given input-output data. In this setting, the contributions of external perturbations and the natural drift of the system to its dynamic behavior cannot be disambiguated, which is critical in many applications. 

In this work, we consider the problem of delineating the contribution of the drift and external perturbations to the dynamic evolution for a class of nonlinear systems. In particular, we consider systems with input-affine dynamics and propose a perturbation strategy to infer the contribution of the control vector field without any influence of the drift dynamics. We show that the resulting learning problem yields a regression function that is an explicit function of the perturbation inputs and therefore directly allows to design experiments to ensure that the regression matrix is well-conditioned. As a result, the proposed perturbation and learning approach, subject to certain requirements on the experimental protocol, can be extremely effective in decoupling the contributions of the external control and the natural drift of a system to its dynamics.

The rest of this paper is organized as follows. In Section \ref{sec:prelims}, we introduce the class of systems considered in this study, and provide a brief background on the model learning problem. We present our approach for data generation and model-learning for recovering the control vector fields independent of the drift vector field in Section \ref{sec:nonlinear_dyn}. We provide several examples and simulation analysis in Section \ref{sec:examples} to demonstrate the effectiveness of the proposed approach. 

\section{System and Data: Model-learning}\label{sec:prelims}
Consider a nonlinear dynamical system with the external control input linearly entering the dynamic equation as described by 
\begin{align}\label{eq:affine_dyn}
    \dot{x}(t) =f(x)+\sum_{k=1}^m g_k(x)u_k(t), \quad x(t_0) \in \Omega,
\end{align}
where $x=(x_1,\ldots,x_n)'\in \Omega \subseteq \mathbb{R}^n$, $t\in [t_0,t_f]$ for $0\le t_0<t_f<\infty$, $u=(u_1,\ldots,u_m)'$, $u_k\in \mathbb{R}$ for each $k \in \{1,\ldots,m\}$ and $f: \Omega \to \mathbb{R}^{n}$ ,  $g_k: \Omega \to \mathbb{R}^{n}$, for $k=1,\ldots,m$, are the drift- and the control- vector field, respectively. We assume that the state evolving on the $n$-dimensional manifold $\Omega$ is accessible for measurement and $f,g_k$ are continuous functions. We will denote
\begin{align}
	\label{eq:g}
	g(x)=(g_1(x),\cdots,g_m(x)) = \Bigg(\begin{smallmatrix}g_{11}(x) & \cdots & g_{1m}(x) \\ \vdots & \ddots & \vdots \\ g_{n1}(x) & \cdots & g_{nm}(x) \end{smallmatrix}\Bigg).
\end{align}

In the rest of this section, we will present a brief overview of the model-learning problem. 
\subsection{Model-learning: A brief overview} 
The primary objective in a system identification or model-learning framework is to reconstruct the dynamics of the system from the given time-series data. For example, given the trajectories $x(t)$ and $u(t)$ the goal is to infer the unknown functions $f,g$ in \eqref{eq:affine_dyn}. This problem has been an active research topic and a recent survey on this topic can be found in \cite{schoukens2019nonlinear}. 

In this section, we will confine our overview to the traditional system identification approaches and some popular learning approaches such as the dynamic mode decomposition \cite{proctor_16_dmdc}. As detailed in \cite{schoukens2019nonlinear}, the existing framework to build dynamical models using the measured input-output data are composed of three major steps: the data collection, the selection of a set of candidate models, and the learning/estimation algorithm to identify the model that best fits the given data. In the following, we illustrate these ideas, and motivate the problem of decoupling control coefficients from the drift dynamics in a model learning framework and point-out the shortcomings of the existing approaches in dealing with this problem.  

\paragraph{Data-generation} 
To infer the dynamics of the system \eqref{eq:affine_dyn}, i.e., unknown functions $f$ and $g$ from data, first, experiments are performed to assimilate input-output data $x(t)$ and $u(t)$. In particular, several samples of data are generated by applying various perturbation inputs, for example, pulse, multi-sine, random noise, etc. \cite{schoukens2019nonlinear}, to the system. The data generated through experiments are fundamental source of information regarding the underlying system dynamics, and diverse data trajectories are desired for learning a reliable model. In practice, the type of excitation signals applied as inputs during data-generation, the sampling rate of the sensors and measurement noise have a major impact on the `quality' of the data obtained, and there is no standard excitation signal that would work for every application \cite{schoukens2019nonlinear}.  

\paragraph{Candidate models- Linear-in-parameter (LIP) form}
The data generated is then used to obtain an approximation of the dynamics from candidate models. A common approach to do this task is to derive an LIP representation of the unknown dynamics, i.e., the system dynamics are re-written in a parametric form as 
\begin{align}\label{eq:lip}
    \dot{x}(t) \approx \Theta \Phi(x(t),u(t)),
\end{align}
where $\Theta = [\Theta_x \quad \Theta_u]$ constitutes the unknown coefficients corresponding to the regression function $\O = [\Phi_x \quad \O_u]'$, each with appropriate dimensions. The input-output data ($x,u$) are then used in \eqref{eq:lip} to estimate the coefficients $\T$ in a supervised learning framework. The choice of regression function $\O$ in \eqref{eq:lip} is a design choice and is typically composed of basis or kernel functions or time-delayed states and control in the discrete-time setting \cite{schoukens2019nonlinear}. 

\paragraph*{Linear dynamical system}
Consider a linear dynamical system with the governing equation given by
\begin{align}\label{eq:linear_control_coefficient}
    \dot{x}(t) = A x(t)+\sum_{k=1}^m b_ku_k(t), \quad x(t_0) \in \mathbb{R}^n,
\end{align}
where $x \in \mathbb{R}^n$, $A \in \mathbb{R}^{n \times n}$ is a linear time-invariant matix, and $b_k \in \mathbb{R}^{n}$ for $k=1,\ldots,m$ is a constant control-coefficient vector. We denote $B = (b_{ij}) =[b_1, \ldots, b_m]\in \mathbb{R}^{n\times m}$ with $i=1,\ldots,n$ and $j=1,\ldots,m$.

For the linear models of the form \eqref{eq:linear_control_coefficient}, one can expect the trajectories $x(t)$ to lie on a hyper-plane. Therefore, a state-space structure of $n$ dimensions can be chosen in a parametric form as in \eqref{eq:lip}, which can then be updated to fit the data samples. For instance, the matrices  $\Theta_x$ and $\Theta_u$ are the approximates of $A$ and $B$, respectively, and the regression matrices $\O_x, \O_u$ correspond to the data samples $x$ and $u$, respectively. The slope of the states $x$ is approximated as $\dot{x}(t) \approx \frac{x(t+t_s)-x(t)}{t_s}$, where $t_s$ is the sampling time, and a linear regression problem is solved to determine the entities of the matrix $\T$. For $t_s$ close to $0$ and sufficiently large data set, the $n^2+nm$ unknown parameters in $\T$ are estimated through solving an associated linear regression problem (LRP).

However, the process of choosing the regression function $\O$ in \eqref{eq:lip} to approximate a nonlinear dynamical system can be quite tedious as illustrated by the following example.

\emph{Example 1:}  A neuron can be modeled using a bio-physically meaningful Hodgkin-Huxley (HH) model, which is a $4$-dimensional model describing the evolution of transmembrane potential in the neuron and it is given by
\begin{align}\label{SI equ: Hodgkin-Huxley}
	C_m\dot{V}_m(t) = &-(\bar{g}_K n_1^4 (V_m-V_K) + \bar{g}_{N_a} n_2^3 n_3(V_m - V_{N_a}) \nonumber \\ & + \bar{g}_l(V_m -V_l)) + I(t),  
\end{align}
and the internal states are governed by 
	\begin{align}
	&\dot{n}_i(t) = \alpha_{n_i}(V_m) (1-n_i(t)) - \beta_{n_i}(V_m)n_i(t), \quad i=1,2,3, \nonumber
 	\end{align}
where $V_m, n_i,\in \mathbb{R}$ are the state variables, $I$ is the control input (extrinsic current, e.g., optogenetic stimulation \cite{williams2015optogenetic}) and $\alpha_{n_i}, \beta_{n_i}$ are nonlinear functions \cite{izhikevich_2007_dynamical}. In many applications of brain medicine, the parameters that model the ionic channel conductances $\bar{g}_k, \bar{g}_{Na}, \bar{g}_l$, the reversal potentials $V_{Na}, V_K, V_l$ and the membrane capacitance $C_m$ are required to be known accurately for the purpose of diagnosis and for synthesis of extrinsic stimulation to evoke specific patterns such as synchronization/desynchronization, etc. \cite{izhikevich_2007_dynamical}. 

In practice, only $V_m$ can be measured, and in this case, the voltage dynamics for the HH model can be rewritten in a parametric form as in \eqref{eq:lip} with the regression function
\begin{align} {\O}(t)=(-n_1^4V_m,n_1^4,-n_2^3n_3V_m,n_2^3n_3,-V_m,1,I)',
\label{eq:example_reg} 
\end{align} 
and the unknown parameters 
\begin{align} \Theta=\frac{1}{C_m}({\bar{g}_{k}},{V_{K}},{\bar{g}_{Na}},{V_{N_a}},{\bar{g}_{l}},{V_l},1).
\label{eq:example_par} \end{align} 
In this example, the selection of candidate models can be very challenging since the gating variables $n_i$ cannot be measured, and even for a fixed set of gating variables $n_i$, the excitation signal $I(t)$ should be designed such that the resulting regression function \eqref{eq:example_reg} is well-conditioned, which is a challenge due to the leak channel and the rich dynamic behavior of a neuron (periodic spiking, bursting, slow-adaptation and fast spiking, etc.) \cite{JCNS_2019biophysically,izhikevich_2007_dynamical}.

\paragraph*{Nonlinear paradigm} For a nonlinear system, the regression function $\O$ in \eqref{eq:lip} is typically composed of orthonormal basis functions of the underlying function space. For instance, if there is a prior knowledge that the system exhibits oscillatory behavior, the Fourier basis can be used to approximate the nonlinear functions. Using a truncated Fourier basis expansion, the unknown dynamics is still represented in the LIP form and the data is used in a supervised learning framework to infer the coefficients of these Fourier terms. Alternatively, polynomials such as Legendre and Chebyshev polynomials are also common choices of candidate regression functions \cite{schoukens2019nonlinear}. 

To enhance the tractability in model selection, it is critical to specify the intended use of the model before experimental design to collect data. For instance, to obtain a predictive model for the purpose of control synthesis \cite{schoukens2019nonlinear}, interpretability of the approximation may not be essential. In this case, abstract nonlinear models such as artificial neural networks or wavelets \cite{ks1990identification,ljung1999system} can be used to approximate the system dynamics in contrast to the LIP-based approximators as in \eqref{eq:lip}. On the other hand, to obtain a simulation model with properties of interpretability, it is essential to carefully select candidate models and design experiments to collect necessary data.

\paragraph{Estimation and learning} Finally, with the parameterization as in \eqref{eq:lip}, the unknown parameters $\T$ are updated to fit the data. In particular, this leads to a linear regression problem. The solution to this problem relies on the choice of regression function $\O$, the complexity of the dynamics of the system and the given data samples. Many variations to this generic model learning framework have been proposed in the literature and a detailed review of these approaches can be found in \cite{schoukens2019nonlinear}.


Due to the LIP representation of the model as in \eqref{eq:lip} and the regression functions $\Phi_x$ and $\Phi_u$, which are typically state-dependent, the contributions of the natural drift of the system and the external input scaled by the control coefficient cannot be tractably decoupled. In the following, we provide some examples to point out the importance of decoupling the drift- and control- vector fields in a model learning problem.

\paragraph*{Example 2} Many control systems describing periodic activity or oscillatory behavior, for example, spiking activity of neurons and thermostatically coupled loads \cite{izhikevich_2007_dynamical, bomela2018phase,kuritz2018ensemble}, can be represented using the phase-reduction model of the form,
\begin{align}\label{eq:phase model}
    \dot{\theta}(t) = \omega(t)+g(\theta(t))u(t)+\eta(t),
\end{align}
where $\t$ is the phase of the oscillations, $\omega$ is the frequency of oscillation, $\eta$ models the (Gaussian) noise and $g$ describes the phase-response curve (PRC). In order to efficiently design control strategies for steering a higher dimensional oscillatory system using its first-order phase model, it is essential to precisely infer the PRC in the model that describes the infinitesimal effect of external weak forcing on the oscillatory behavior. In practice, the time-varying drift and the noise can affect the computation of the PRC, and it is essential to delineate the contribution of the time-varying frequency from the PRC in order to design precise control signals.

Furthermore, in many practical applications, actuators contribute non-trivially to the dynamics of the system. For instance, for the HH model \eqref{SI equ: Hodgkin-Huxley} in \emph{Example 1}, the perturbation $I(t)$ is typically applied using external apparatus such as optogenetic stimulation \cite{grossman2011modeling,williams2015optogenetic}, which introduces nonlinearity and the control vector field is no longer a constant determined by the membrane capacitance ($C_m$) as in \eqref{SI equ: Hodgkin-Huxley}. Therefore, it is imperative to disambiguate the contribution of the input channel from the drift (ionic channels) to the overall dynamics of a neuron in order to enable interpretable inference of the underlying dynamical process.  

In view of the above practical challenges arising in emerging applications, we consider the problem of disambiguating the drift and control vector fields in a model inference problem. In particular, we propose a methodology to exploit the structure of input-affine nonlinear system to derive a regression problem for inferring the control vector field, which is independent of the drift vector field. To the best of our knowledge, such a problem is not considered in the existing model learning/system identification algorithms.   

In the next section, we will consider the nonlinear input-affine system and present our formulation for learning the control vector field independent of the drift vector field.
\section{Nonlinear dynamical systems}\label{sec:nonlinear_dyn}
We propose a data generation and learning procedure to recover the control coefficient $g$. In particular, we propose to perform multiple experiments under different experimental settings to generate data samples over a short time horizon and formulate a learning problem for reconstruction of $g$ such that the learning problem does not depend on $f$. Consider the input-affine nonlinear system given by \eqref{eq:affine_dyn}, we assume that $x$ is accessible for measurements and the system can be perturbed using the external control input $u(t)$. 

\paragraph*{Perturbation strategy} 

Let the $k^{th}$ perturbation input to the $i^{th}$ experiment be denoted by $u_k^{(i)}(t)$ for $i=0,1,2,\ldots,{N}$ and $k=1,2,\ldots,m$. 
Following this perturbation, 
the resulting system dynamics are described by
\begin{align}\label{eq:ith_nonlinear_control_coeff}
    &\dot{x}^{(i)}(t) =  f(x(t))+\sum_{k=1}^m g_k(x(t)) u_k^{(i)}(t),\ \nonumber \\ & x^{(i)}(t_0)=a_0, \quad i=0,\ldots,N.
\end{align}

Let the difference between the applied control in the $0^{th}$ (or the reference) and the $i^{th}$ experiment be defined as $\Delta_{i} u_k(t) = u_k^{(0)}(t)-u_k^{(i)}(t)$, and let the difference in the resulting state dynamics be denoted as $\Delta_{i} \dot{x}(t) =  \dot{x}^{(0)}(t) -  \dot{x}^{(i)}(t), \quad i=1,\ldots,N$, where $x^{(i)}(t_0)=a_0$ in every experiment. By setting the same initial condition in each experiment (i.e., $t_0$ and $x^{(i)}(t_0)$ for $i=0,1\ldots,N$), we have 
\begin{align}\label{eq:g_regression_prelim}
&\Delta_{i} \dot{x}(t_0) = g(x(t_0)) \Delta_{i} u(t_0), \quad \Delta_{i} {x}(t_0)=0,
\end{align}
where $g$ is defined in \eqref{eq:g}. The equation in \eqref{eq:g_regression_prelim} reveals a system of linear equations, which can be exploited to determine the nonlinear function $g$ at $a_0$. 

To approximate $g$ in its domain $\Omega$, we repeat this experimental procedure at different sample points in $\Omega$. In particular, let $a_1,a_2,\ldots,a_{M}$ be distinct points in $\Omega$. We apply $u_k^{(i)}$ for $k=1,\ldots,m$ and $i=0,\ldots,N$ at each of the initial states $a_1,\ldots,a_{M}$ resulting in a total of $(M+1)(N+1)$ experiments. The dynamics of the system under the proposed experimental setting can be represented using \eqref{eq:ith_nonlinear_control_coeff} as
\begin{align}\label{eq:ith_nonlinear_control_coeff_gener}
    & \dot{x}^{(i)}(t) = f(x(t))+\sum_{k=1}^m g_k(x(t)) u_k^{(i)}(t),\nonumber \\ & x^{(i)}(t_0)=a_j, \quad i=0,\ldots,N, \quad j=0,\ldots,M.
\end{align}
Without loss of generality, set $t_0=0$ and define $\Delta_{i} \dot{x}(0,a_j) = (\Delta_{i} \dot{x}_1(0,a_j), \ldots, \Delta_{i} \dot{x}_n(0,a_j))' $ to denote the derivative at time $t=0$ with the initial condition $a_j$, and denote the $m$-dimensional control at $t=0$ as $\Delta_{i} u(0) = (\Delta_{i} u_1(0),\ldots, \Delta_{i} u_m(0))'$. This leads to a system of linear equations of the form 
\begin{align}\label{eq:lips_control_affine}
&\Delta_{i} \dot{x}_j(0,a_k) = \sum_{s = 1}^m g_{js}(a_{k}) \Delta_{i} u_s(0),  \\ & i=1,\ldots, N, \quad j=1,\ldots,n, \quad k=0,\ldots, M.\nonumber
\end{align}
where $g_{js}$ is defined in \eqref{eq:g}. Using the data generated by the proposed strategy, we can reconstruct or learn an approximation of the function $g$ from \eqref{eq:lips_control_affine} using truncated orthonormal bases, such as the Legendre polynomials or the Fourier basis, i.e.,
{$g_{ij}(x) \approx \sum_{k=0}^{L_{ij}} \alpha_{ij}^{(k)} \phi_{ij}^{(k)}(x)$, $i = 1, \ldots, n$, $j = 1, \ldots, m$,}
where $\{\a_{ij}^{(k)}\}_{k=0}^{L_{ij}}$ are scalar coefficients, $\{\phi_{ij}^{(k)}\}_{k=0}^\infty$ are orthonormal basis functions and $L_{ij}$ are the number of expansion terms in the truncated series $\{\phi_{ij}^{(k)}\}_{k=0}^{L_{ij}}$ for $i=1,\ldots,n$ and $j=1,\ldots,m$. Such an approximation is possible as any continuous function can be approximated arbitrarily well on a compact support using orthonormal bases by the Stone-Weierstrass theorem \cite{rudin1976principles}. 
This leads to 
\begin{align}\label{eq:n_lips_affine}
&\Delta_{i} \dot{x}_j(0,a_p) \approx \sum_{s=1}^m \sum_{k=0}^{L_{js}} \alpha_{js}^{(k)} \phi_{js}^{(k)}(a_p)\Delta_i u_s(0),\nonumber\\
& i=1,\ldots, N, \quad p=0,\ldots, M, \quad j=1,\ldots,n.
\end{align}


\begin{remark}
The $n$ linear regression problems (LRPs) resulting from \eqref{eq:n_lips_affine} can be solved (in parallel) for the unknown coefficients $\a_{js}$ to obtain an approximation for $g$. Note here that the solvability of these LRPs is directly related to the variation in the perturbation inputs applied during each experiment, choice of basis functions and the short-term evolution of the system, and it is independent of the natural drift $f$ of the system. 
\end{remark}
\subsection{Revisiting the case of linear dynamical systems}\label{sec:linear_dyn}
Here, we will revisit the case when the dynamics of the system are linear as in \eqref{eq:linear_control_coefficient}. Suppose that the matrix $A$ and the vectors $b_k$ are unknown and that the state $x$ is accessible for measurement. 

To recover the control-coefficient $B$, we need to perform $N+1$ experiments wherein we apply different perturbation inputs $\{u_k^{(0)}(t),\ldots,u_k^{(N)}(t)\}$ for $k=1,\ldots,m$, and record the resulting states. We can then analyze the dynamical equation \eqref{eq:linear_control_coefficient} for different control signals. In particular, the system dynamics driven by control inputs in the $i^{th}$ experiment is given by
\begin{align}
    \dot{x}^{(i)}(t) = & Ax(t)+\sum_{k=1}^m b_k u_k^{(i)}(t), \quad x(t_0)=x_0.
\end{align}
We also know the difference between the applied control in the reference ($0^{th}$) experiment and the $i^{th}$ experiment,
$\Delta_{i} u_k(t) = u_k^{(0)}(t)-u_k^{(i)}(t)$,
which results in $\Delta_{i} \dot{x}(t) =  \dot{x}^{(0)}(t) -  \dot{x}^{(i)}(t)$, for $i=1,\ldots,N$. 
Unlike the recovery of nonlinear $g$ in \eqref{eq:ith_nonlinear_control_coeff_gener},  for the linear time-invariant control coefficients, we have
\begin{align}\label{eq:lip_linear_cont_coeff}
\Delta_{i} \dot{x}(t_0) &= \sum_{k=1}^m b_k \Delta_{i} u_k(t_0) ,  \\ \quad \Delta_{i} {x}(t_0)&=0, \quad i=1,\ldots,N.\nonumber
\end{align}
Without loss of generality and for ease of exposition, we set $t_0=0$, and define $\Delta_{i} \dot{x}(0) = (\Delta_{i} \dot{x}_1(0), \ldots, \Delta_{i} \dot{x}_n(0))' $ and $\Delta_{i} u(0) = (\Delta_{i} u_1(0),\ldots, \Delta_{i} u_m(0))'$. Using these definitions in \eqref{eq:lip_linear_cont_coeff} leads to $n$ LRPs, given by

\begin{align}\label{eq:jth_lip_linear_cont_coeff}
    \left[ \begin{smallmatrix}\Delta_{1}u_1&\cdots&\Delta_{1} u_m \\
                  \Delta_{2}u_1&\cdots&\Delta_{2} u_m \\
                  \vdots &\ddots &\vdots \\ 
                  \Delta_{N} u_1& \cdots & \Delta_{N} u_m\end{smallmatrix}\right] \left[\begin{smallmatrix}b_{j1} \\ b_{j2} \\ \vdots \\ b_{jm} \end{smallmatrix} \right] = \left[\begin{smallmatrix}\Delta_{1} \dot{x}_{j} \\ \Delta_{2} \dot{x}_{j} \\ \vdots \\ \Delta_{N} \dot{x}_{j} \end{smallmatrix} \right],
\end{align}
for $j=1,\ldots, n$ at $t=t_0 (=0)$. In particular, for $N=m$, the perturbation inputs can be explicitly designed to ensure that regression matrix in \eqref{eq:jth_lip_linear_cont_coeff} is always full-rank, which will yield an exact solution to the associated LRPs. For numerical computation, the derivatives in \eqref{eq:jth_lip_linear_cont_coeff} can be computed by using the forward difference, $\dot{x}_{j}(t_0) \approx \frac{{x}_{j}(t_0+t_s)-{x}_{j}(t_0)}{t_s}$ with $t_s$ being the sampling time, which inherits numerical errors on the order of $\mathcal{O}(t_s)$.
\begin{remark}
Once the control vector field is inferred using the proposed method, the regression problem for inferring the drift-vector field can be formulated in a supervised learning framework of the form $Ax = b$, where the output $b = \dot{x}(t)-\sum_{i=1}^m g_i(x)u_i(t)$, $A$ is composed of the basis functions evaluated at the sample points, and $x$ constitutes the unknown coefficients. Note that in the presence of unmodeled fluctuations or noise in the dynamics, the linear equations in \eqref{eq:lip_linear_cont_coeff} (or in \eqref{eq:lips_control_affine}) will be of the form $Ax +\d = b$, where $\d$ models the fluctuations. Under the assumption that these fluctuations are modeled by white noise, the least-square solution to the regression problem will be the maximum-likelihood estimate and as the number of experiments increase, the solution will converge in probability.
\end{remark}

In the next section, we apply the proposed method to recover the control vector fields for linear, bilinear, and nonlinear dynamical systems.

\section{Simulation Results}\label{sec:examples}
In this section, we demonstrate the effectiveness of the proposed approach using three numerical examples.
\subsection{Linear dynamical systems}
Here, we consider a second-order linear dynamical system to analyze the proposed approach and to illustrate its advantages. In particular, we study one taking the form in \eqref{eq:linear_control_coefficient} with the matrices $A = \left[\begin{smallmatrix}1&4 \\ 5&-1\end{smallmatrix}\right]$ and $B = \left[\begin{smallmatrix}2&1 \\ 0.6&1\end{smallmatrix}\right]$. We fixed the initial state $x(t_0)=(0,-0.25)'$ at $t_0=0$ and performed $3$ experiments with respective controls $u^{(i)}(t_0)=(1,2)',(2,4)',(3,8)'$ for $i=0,1,2$ and $t_s=1 ms$.
\begin{figure}
    \centering
    \includegraphics[width=0.93\linewidth,keepaspectratio]{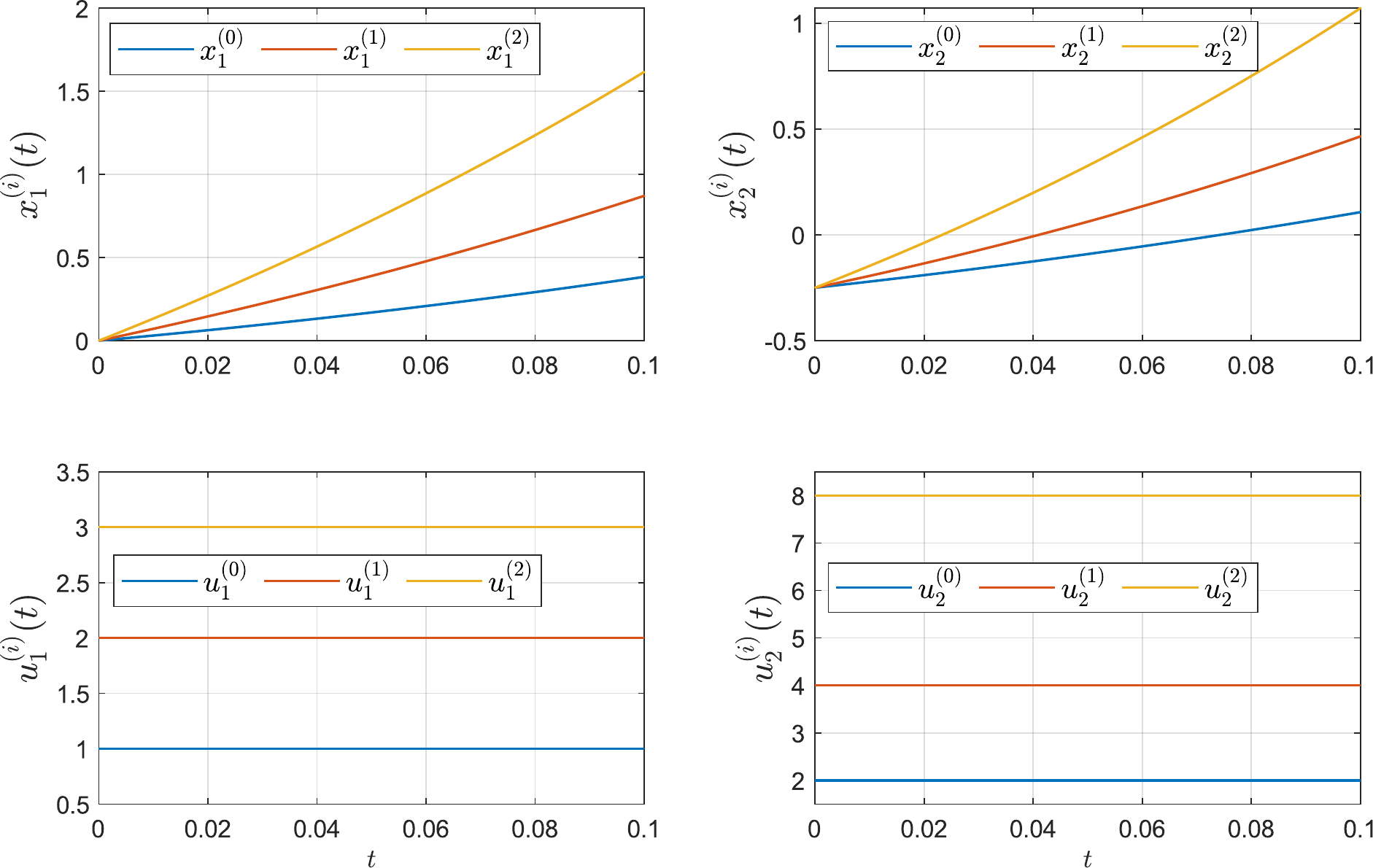}
    \caption{The state and control trajectories recorded during the $3$ data generating experiments. The initial time $t_0=0$ and the initial states for all the three experiments are $(0,-0.25)'$.}
    \label{fig:lin_experiments}
\end{figure}
This resulted in the following equations,
\begin{align}
    \dot{x}^{(0)}(0) = & \big[\begin{smallmatrix} -1 \\ 0.25 \end{smallmatrix}\big]+\big[\begin{smallmatrix} 4 \\ 2.6 \end{smallmatrix}\big], \quad
    \dot{x}^{(1)}(0) = & \big[\begin{smallmatrix} -1 \\ 0.25 \end{smallmatrix}\big]+\big[\begin{smallmatrix} 8 \\ 5.2 \end{smallmatrix}\big], \nonumber \\
    \dot{x}^{(2)}(0) = & \big[\begin{smallmatrix} -1 \\ 0.25 \end{smallmatrix}\big]+\big[\begin{smallmatrix} 14 \\ 9.8 \end{smallmatrix}\big], \nonumber
\end{align}
and the variation equations 
\begin{align}\label{eq:linear_ex_lip}
    \Delta U \left[\begin{smallmatrix}b_{11} \\ b_{12} \end{smallmatrix} \right] = \left[\begin{smallmatrix}\Delta_1 \dot{x}_{1}(0) \\ \Delta_2 \dot{x}_{1}(0) \end{smallmatrix} \right], \
    \Delta U \left[\begin{smallmatrix}b_{21} \\ b_{22} \end{smallmatrix} \right] = \left[\begin{smallmatrix}\Delta_1 \dot{x}_{2}(0) \\ \Delta_2 \dot{x}_{2}(0) \end{smallmatrix} \right],
\end{align}
with $\Delta U = \left[ \begin{smallmatrix}\Delta_1u_1(0)&\Delta_1 u_2(0) \\
                  \Delta_2u_1(0)&\Delta_2 u_2(0)\end{smallmatrix}\right] = \left[ \begin{smallmatrix}-1&-2 \\
                  -2&-6\end{smallmatrix}\right]$ and 
                  $$\left[\begin{smallmatrix}\Delta_1 \dot{x}_{1}(0) \\ \Delta_2 \dot{x}_{1}(0) \end{smallmatrix} \right] \approx \left[\begin{smallmatrix}-4.0007 \\ -10.0018 \end{smallmatrix} \right], \left[\begin{smallmatrix}\Delta_1 \dot{x}_{2}(0) \\ \Delta_2 \dot{x}_{2}(0) \end{smallmatrix} \right] \approx \left[\begin{smallmatrix}-2.6009 \\ -7.2021 \end{smallmatrix} \right].$$

The unknown $B$ matrix is then computed from \eqref{eq:linear_ex_lip} as $\left( \begin{smallmatrix}2.0002&1.0003 \\ 0.6005&1.0002\end{smallmatrix}\right)$. An illustration of the state and control trajectories used to recover the $B$ matrix is given in Fig. \ref{fig:lin_experiments}.

\subsection{Bilinear systems} A Bloch system described by a third-order model has its trajectories evolving on a unit sphere with the dynamics given in \eqref{eq:affine_dyn} with $m=2$, 
$$f(x) = \left[\begin{smallmatrix} -\omega x_2\\ \omega x_1\\0 \end{smallmatrix}\right], g_1(x) = \left[\begin{smallmatrix} \varepsilon x_3\\0\\-\varepsilon x_1\end{smallmatrix}\right], g_2(x) = \left[\begin{smallmatrix} 0\\-\varepsilon x_3\\\varepsilon x_2\end{smallmatrix}\right],$$
where $\varepsilon=0.6$ and $\omega =  1.4$. The unforced trajectories of this system will only contain limited information regarding the underlying system structure, which cannot be used to infer a reliable model for the system. Therefore, it is essential to design efficient excitation signals to extract information of the system structure. 

\begin{figure}[h]
    \centering
    \includegraphics[width=\linewidth,keepaspectratio]{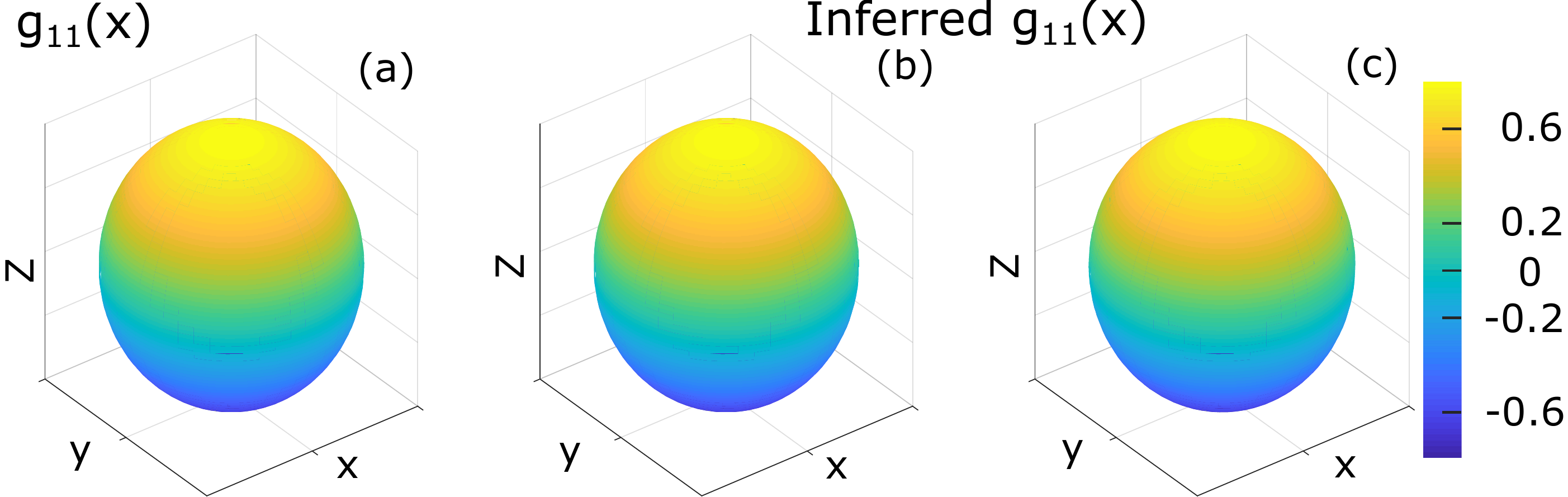}
    \caption{Bloch system: To validate the approximation of $g_1(x),g_2(x)$, we randomly sampled over $1000$ sample points on the sphere that were not used as initial conditions when performing the experiments and evaluated the function $g_1 \in \mathbb{R}^{3}, g_2\in \mathbb{R}^{3}$ at these sample points, and the resulting values of $g_{11}(x)$ are recorded. (a) The actual values of $g_{11}(x)$; (b) Inferred $g_{11}(x)$ with polynomial basis functions; and (c) Inferred $g_{11}(x)$ with Fourier basis functions. The color bar denotes the range of $g_{11}(x)$ in the domain $\Omega$.}
    \label{fig:bloch}
\end{figure}

We consider two cases. First, we define the regression function using the first and the second order polynomials of the states, i.e., $(1,x_1,x_2,x_3,x_1^2,x_2^2,x_3^2,x_1x_2,x_1x_3,x_2x_3)'$. In the second case, we define the regression function using the Fourier basis with $L_{js}=5$ for $j=1,2,3, s=1,2$ in \eqref{eq:n_lips_affine}. We randomly selected $20$ initial states and applied controls of the form $(k,0),(0,k)$ for $k=0,1,\ldots,3$ to recover $g_1$ and $g_2$. The coefficients estimated using the proposed approach resulted in the replicates of the function $g_1,g_2 \in [-1,1]$ as shown in the representative Fig. \ref{fig:bloch}, where the function $g_{11}(x)=\varepsilon x_3$ is recorded, and its approximation using the polynomial basis and the Fourier basis are also recorded. 
\subsection{Nonlinear Oscillators} 
Here, we revisit the phase model given in \emph{Example 3} and apply the proposed approach for recovering the PRC. We consider the dynamics as in \eqref{eq:phase model} with the PRC given as $g(\t) = -\sin(\t)\exp(3[\cos(\t - 0.9\pi) - 1])$, $\omega(t)= 0.1t$ and $\eta=0$.

\begin{figure}[h]
    \centering
    \includegraphics[width=0.9\linewidth,keepaspectratio]{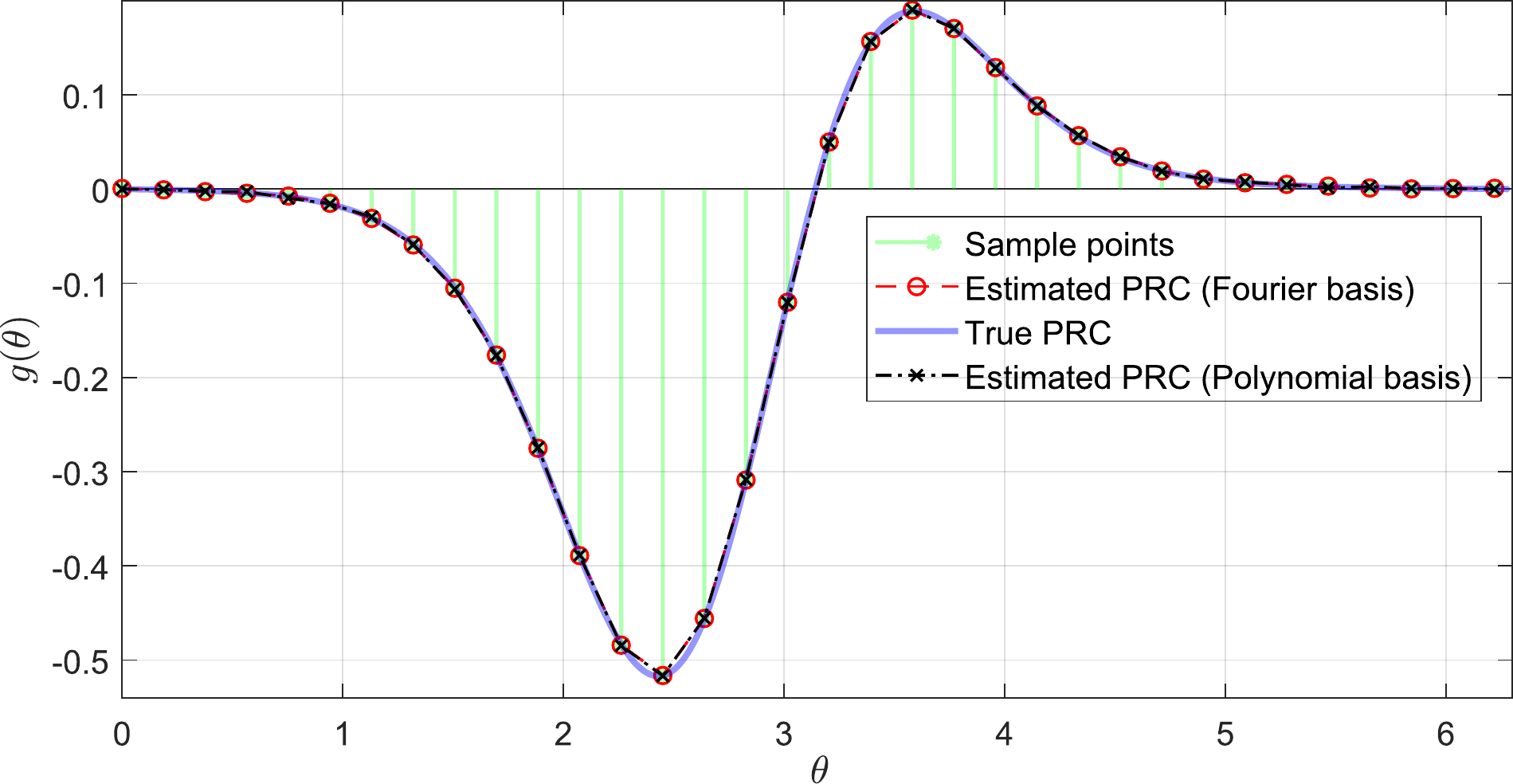}
    \caption{Actual PRC ($g(\t)$) recorded in blue solid line. The inferred PRC with Fourier basis and polynomial basis are recorded with red and black dashed lines, respectively ($u(t)\in [-1,1]$). We performed 3 experiments at 35 sample points, randomly selected from $[0,2\pi]$, and used the resulting data was used to recover the Fourier coefficients.}
    \label{fig:prc}
\end{figure}

\begin{figure}[h]
    \centering
    \includegraphics[width=0.9\linewidth,keepaspectratio]{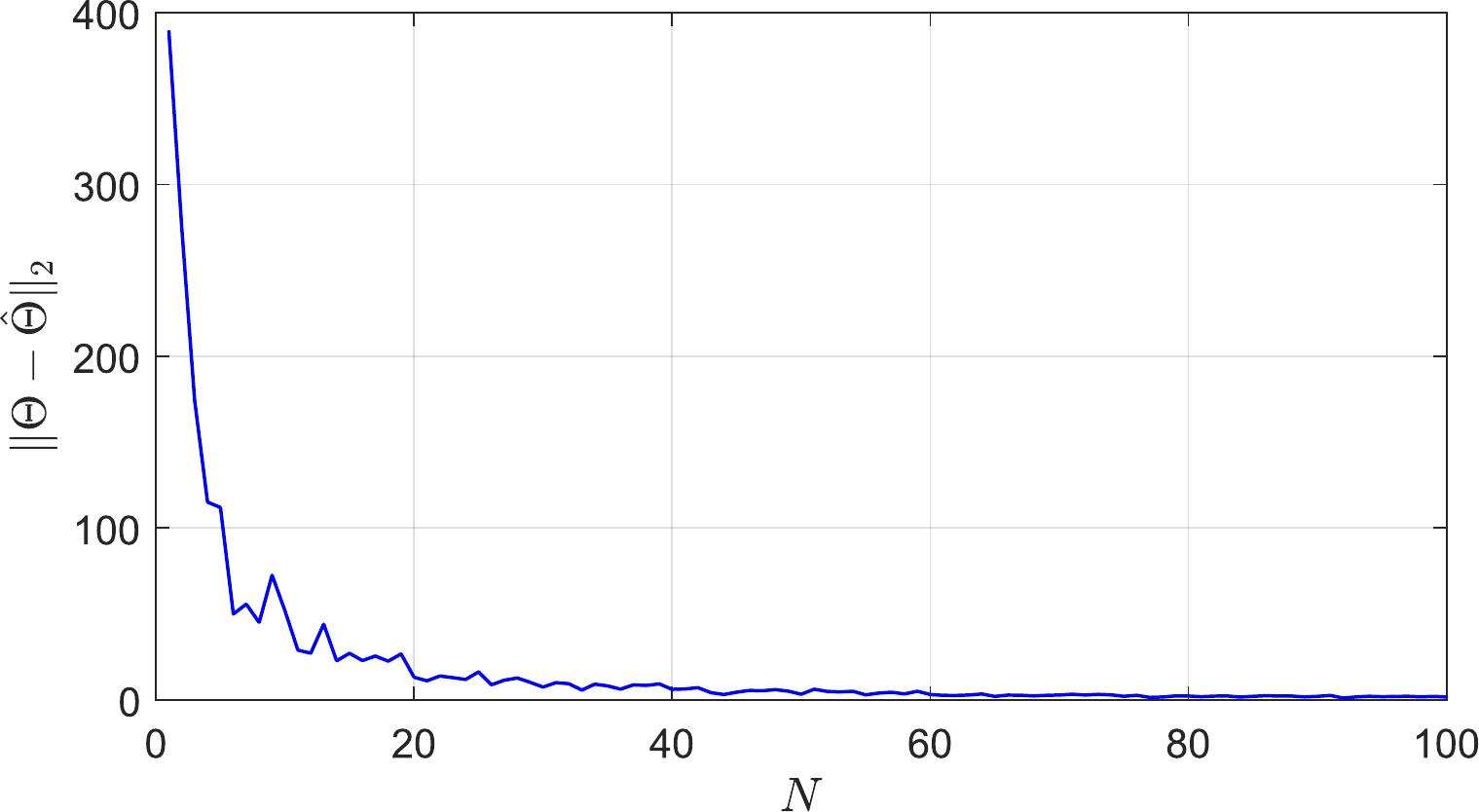}
    \caption{Parameter estimation error vs the number of perturbation experiments. Here $\T$ denotes the Fourier coefficients obtained with noise free data (corresponding to the inferred PRC in Fig. \ref{fig:prc}) and $\hat{\T}$ denotes the estimated Fourier coefficients with internal fluctuations in the dynamics ($\eta \in [-1,1]$). The proposed data-generation and approximation is performed for varying number of perturbation experiments $N$, and the norm of the parameter error is recorded.}
    \label{fig:prc_error}
\end{figure}

We use the truncated Fourier basis expansion to approximate the PRC upto order $6$. The actual and recovered PRCs and the sample points (total sample points $35$) are recorded in Fig. \ref{fig:prc}. Additionally, we considered the case when the dynamics includes an unmodeled white noise as in \eqref{eq:phase model} ($\eta \in [-1,1]$). In this case, as the number of perturbation experiments increased, the parameter estimation error, which was computed for the purpose of illustration using the difference between the Fourier coefficients obtained in the noise free case and the Fourier coefficients estimated in the presence of noise, converged in to zero as shown in Fig. \ref{fig:prc_error}.

\section{Conclusions}\label{sec:conclusions}
In this work, we have proposed a data generation protocol and learning framework for recovering the control vector fields, wherein the learning problem of recovering the control vector fields does not depend on the natural drift of the system. Using the input-affine structure, we proposed a perturbation strategy and demonstrated that the control and drift vector fields can be decoupled in a data-driven modeling framework, which under certain conditions pertaining to availability of access to actuation of the system at arbitrary initial conditions, can yield excellent performance when enacted on unknown nonlinear dynamical systems. We demonstrated the efficacy of our approach by implementing the proposed methods in several numerical examples.
\bibliographystyle{IEEEtran}
\bibliography{CDC_20.bib}
\addtolength{\textheight}{-12cm}   
\end{document}